\documentclass[12pt]{article}
\usepackage{amsfonts,amsmath,amssymb,setspace,mathtools,bm}
\usepackage{graphicx,float,cleveref}
\usepackage{tikz}

\newtheorem{Theorem}{Theorem}
\newtheorem{Lemma}[Theorem]{Lemma}
\newtheorem{Corollary}[Theorem]{Corollary}

\newenvironment{Proof}{\begin{trivlist} \item[] {\bf Proof.}}{\hfill $\Box$\end{trivlist}}

\renewcommand{\geq}{\geqslant}
\renewcommand{\leq}{\leqslant}

\title{Existential Closure in Line Graphs}

\author{Andrea C. Burgess\thanks{Department of Mathematics and Statistics, University of New Brunswick, Saint John, NB, E2L 4L5, Canada, andrea.burgess@unb.ca}
       \and Robert D. Luther\thanks{Department of Mathematics and Statistics, Memorial University of Newfoundland, St. John's, NL, A1C 5S7, Canada, rdl863@mun.ca (corresponding author)}
       \and David A. Pike\thanks{Department of Mathematics and Statistics, Memorial University of Newfoundland, St. John's, NL, A1C 5S7, Canada, dapike@mun.ca}}

\begin{document}

\maketitle

\begin{abstract}
A graph is {\it $n$-existentially closed} if, for all disjoint sets of vertices $A$ and $B$ with $|A\cup B|=n$, there is a vertex $z$ not in $A\cup B$ adjacent to each vertex of $A$ and to no vertex of $B$.

In this paper, we investigate $n$-existentially closed line graphs. 
In particular, we present necessary conditions for the existence of such graphs as well as constructions for finding infinite families of such graphs. 
We also prove that there are exactly five $2$-existentially closed planar line graphs.
We then consider the existential closure of the line graphs of hypergraphs and present constructions for $2$-existentially closed line graphs of hypergraphs.
\end{abstract}


MSC Classification Codes: 05C10, 05C65, 05C76, 05C99

\pagebreak

\section{Introduction}

For a positive integer $n$, a graph with at least $n$ vertices is {\it $n$-existentially closed} or simply {\it $n$-e.c.\ }if, for all disjoint sets of vertices $A$ and $B$ with $|A\cup B|=n$, there is a vertex $z$ not in $A\cup B$ adjacent to each vertex of $A$ and no vertex of $B$. 
Hence, for all $n$-subsets $S$ of vertices, there exist $2^n$ vertices joined to $S$ in all possible ways. 
For example, a 1-e.c.\ graph is one with neither isolated nor universal vertices.

If a graph has the $n$-e.c.\ property, then it possesses other structural properties such as the following.

\begin{Theorem}\cite{Bon}\label{1}
Let $G$ be an $n$-e.c.\ graph where $n$ is a positive integer.

\begin{enumerate}
\item The graph $G$ is $m$-e.c.\ for all $1\leq m\leq n-1$.
\item The graph $G$ has order at least $n+2^n$, and has at least $n2^{n-1}$ edges.
\item The complement graph is $n$-e.c.
\item Each graph of order at most $n+1$ embeds in $G$.
\item If $n>1$, then for each vertex $x$ of $G$, each of the graphs $G-x$, the subgraph induced by the neighbourhood $N(x)$, and the subgraph induced by $(V(G)\setminus N(x))-x$ are $(n-1)$-e.c.
\end{enumerate}
\end{Theorem}

There are exactly three non-isomorphic 1-e.c.\ graphs of order 4: $2K_2$, $C_4$, and $P_4$. 
The cartesian product $K_3\square K_3$ is the unique 2-e.c.\ graph of minimum order \cite{BonCam}. 

When studying graph properties which represent structure between vertices, it is often interesting to consider how similar structure can be described in terms of the edges of a graph. 
Furthermore, when investigating the structure between the edges of a graph, it is often convenient to examine the corresponding line graph.
For a graph $G$, the {\it line graph} of $G$, denoted $L(G)$, is the graph whose vertices are the edges of $G$ and two vertices in $L(G)$ are adjacent if and only if they correspond to adjacent edges in $G$, where two edges are adjacent if and only if they share at least one end vertex. 

One such property with interesting behaviour involves graph colourings.
For example, the {\it chromatic number} of a graph $G$, denoted $\chi(G)$, is the smallest number of colours required to colour the vertices of $G$ such that no two vertices of the same colour are adjacent, whereas the {\it chromatic index} of $G$, denoted $\chi'(G)$, is the smallest number of colours required to colour the edges of $G$ such that no two edges of the same colour are incident with the same vertex.
Note that by definition, for a graph $G$, the chromatic index $\chi'(G)$ is equal to the chromatic number of its line graph, $\chi(L(G))$.

Another structural property with similar behaviour is that of independence.
Let $\alpha(G)$ denote the maximum number of pairwise non-adjacent vertices in a graph $G$, and similarly let $\alpha'(G)$ denote the maximum number of pairwise non-adjacent edges (i.e., $\alpha'(G)$ denotes the size of a maximum matching in $G$). 
It is easy to observe that $\alpha'(G) = \alpha(L(G))$ for any graph $G$.

In a similar spirit, we define a version of existential closure expressed in terms of the edges of a graph.
For a positive integer $n$, we say that a graph is {\it $n$-line existentially closed} or simply {\it $n$-line e.c.}\ if, for all disjoint sets of edges $A$ and $B$ with $|A\cup B|=n$, there is an edge $e$ not in $A\cup B$ adjacent to each edge of $A$ and no edge of $B$. 
By definition, a graph is $n$-line e.c.\ if and only if its line graph is $n$-e.c.\ and so, our main topic will be investigating the $n$-e.c.\ property in line graphs. 
In \cite{HPS}, the notation $\Xi(G)$ was first introduced to represent the largest integer $n$ for which the graph $G$ is $n$-e.c.
We may similarly define $\Xi'(G)$ to be the largest integer $n$ for which the graph $G$ is $n$-line e.c.\ and observe that $\Xi'(G)=\Xi(L(G))$.

The only graphs which fail to be 1-line e.c.\ are those containing an edge that is adjacent to every other edge and disconnected graphs containing a connected component consisting of a single edge.
For this reason, we will only consider connected graphs going forward.
Note that if a graph $G$ has a duplicate edge, then it would be impossible to find a third edge adjacent to one but not the other and therefore $G$ would not satisfy the 2-line e.c.\ property. 
Similarly, if $G$ has a loop at vertex $v$, consider any other edge incident with $v$ and note that no third edge would be adjacent to the loop but not the other edge and so $G$ would not satisfy the 2-line e.c.\ property. 
For these reasons, as we continue our exploration of $n$-line e.c.\ graphs with values of $n$ greater than 1, all graphs we discuss are assumed to be simple.

For small values of $n$, examples of $n$-line e.c.\ graphs are easy to find. 
For instance, the graph $C_4$ is the unique 1-line e.c.\ graph of minimum order. 
The graph $K_{3,3}$ is the unique 2-line e.c.\ graph on 9 edges since it is the only graph with corresponding line graph $K_3\square K_3$. 
Also, complete graphs on at least six vertices and complete bipartite graphs with part sizes at least three are all 2-line e.c.\ graphs.

In this paper, we focus on finding conditions for the existence of $n$-line e.c.\ graphs. 
In Section \ref{Res}, we focus on finding necessary conditions for the existence of such graphs.
In particular, we prove that if $G$ is $n$-line e.c.\ then $n$ is at most 2, narrowing our focus to finding examples of 2-line e.c.\ graphs.
In Section \ref{Cons}, we present constructions which generate infinite families of 2-line e.c.\ graphs and in Section \ref{Planar}, we prove that there are exactly five graphs which are both $2$-line e.c.\ and planar.
Lastly, in Section \ref{hyper}, we introduce the problem of existential closure in the line graphs of hypergraphs and present constructions for 2-line e.c.\ hypergraphs.

\section{$n$-Line e.c.\ Graphs}

When studying a combinatorial object it is often natural to ask what the necessary and sufficient conditions are for such an object to exist.
In this section, we look closely at these conditions in order to ultimately construct more examples of such graphs.

\subsection{Necessary Conditions}\label{Res}

When studying $n$-line e.c.\ graphs, there are some immediate necessary conditions that can be observed from Theorem \ref{1} parts 1 and 2.

\begin{Theorem}\label{2}
If $G$ is an $n$-line e.c.\ graph, then $G$ is $m$-line e.c.\ for all $1\leq m\leq n-1$ and $G$ has at least $n+2^n$ edges.
\end{Theorem}

The following theorem poses a heavy restriction on the existence of $n$-line e.c.\ graphs.

\begin{Theorem}\label{n=2}
Let $G$ be a graph. If $G$ is $n$-line e.c.\ then $n\leq 2$. Alternatively, if $G$ is a graph, then $\Xi'(G)\leq 2$.
\end{Theorem}

\begin{Proof}
By Theorem \ref{2}, we know that an $n$-line e.c.\ graph is also $m$-line e.c.\ for $1\leq m\leq n$. 
So it suffices to show that $G$ cannot be 3-line e.c. 

Suppose $G$ is 3-line e.c.\ and let $e_0$ be an edge of $G$. 
Since $G$ is also 1-line e.c., there exists an edge $e_1$ not adjacent to $e_0$. 
Also, since $G$ is 2-line e.c., there exists an edge $e_2$ not adjacent to $e_0$ or $e_1$. 
Finally, since $G$ is 3-line e.c., there must exist a fourth edge, adjacent to each of the previous three distinct edges, no two of which are adjacent. 
This is impossible, so $G$ cannot be 3-line e.c.
\end{Proof}

As an alternative proof of Theorem \ref{n=2} consider the following: Suppose $G$ is a 3-line e.c.\ graph, that is, $L(G)$ is a 3-e.c.\ graph. 
Note that each graph of order four must occur as an induced subgraph in $L(G)$ by Theorem \ref{1} part 4. 
In particular, $L(G)$ must contain $K_{1,3}$ as an induced subgraph, but this is impossible since all line graphs are necessarily claw-free \cite{Beineke}.

Due to the implication of Theorem \ref{n=2}, our attention will now focus specifically on 2-line e.c.\ graphs. A graph $G$ is 2-line e.c.\ if and only if for each pair of distinct edges $e,f\in E(G)$, the following hold:
\begin{itemize}
\item[$(i)$] there is another edge adjacent to both $e$ and $f$,
\item[$(ii)$] there is another edge adjacent to neither $e$ nor $f$, and
\item[$(iii)$] there is another edge adjacent to $e$ but not to $f$.
\end{itemize}
This simplification of the definition leads to the observation of several consequential properties.

\begin{Lemma}\label{deg3}
Let $G$ be a 2-line e.c.\ graph. Then the minimum degree of $G$, $\delta(G)$, is at least three.
\end{Lemma}

\begin{Proof}
Suppose $x$ is a vertex of degree one in $G$, let $e$ be the edge incident with $x$ and let $f$ be any other edge in the graph. 
Applying condition $(i)$ to $e$ and $f$ gives a third edge $g$ adjacent to $e$ and applying $(iii)$ to $e$ and $g$ forces a second edge to be incident with $x$.

Now suppose $x$ is a vertex of degree two in $G$ and let $e$ and $f$ be distinct edges incident with $x$. 
By applying condition $(i)$ to $e$ and $f$, there must exist a third edge $g$ which forms a triangle with $e$ and $f$.
Now apply condition $(iii)$ to $e$ and $g$ and observe that $x$ must have degree at least three.
\end{Proof}

Furthermore, for graphs with minimum degree at least three, condition $(iii)$ is implied by conditions $(i)$ and $(ii)$.
Suppose $e=\{u,x\}$ and $f=\{v,x\}$ are two adjacent edges in a 2-line e.c.\ graph $G$. 
Since $\delta(G)\geq 3$, vertex $u$ has at least three neighbours at most two of which could be $x$ and $v$.
So there must exist an edge adjacent to $e$ but not to $f$ and vice versa. 
Otherwise, letting $e$ and $f$ be two disjoint edges, condition $(ii)$ implies there exists a third edge $g$ adjacent to neither $e$ nor $f$.
Now apply condition $(i)$ to $e$ and $g$ and observe that this edge is adjacent to $e$ but not to $f$.
Therefore, when checking a graph for the 2-line e.c.\ property, it is sufficient to verify only conditions $(i)$ and $(ii)$.

Additionally, if $G$ is a 2-line e.c.\ graph, then every matching of size two in $G$ can be embedded in a matching of size three (where a matching is a set of disjoint edges), as well as a path of length three.
Consequently, $G$ has no induced matching of size two.
Such graphs are often referred to as $2K_2$-free graphs; for characterisations of these graphs see \cite{DSM} and \cite{Meister}.
Therefore, the class of 2-line e.c.\ graphs is contained within the class of $2K_2$-free graphs.

\subsection{Constructing 2-Line e.c.~Graphs}\label{Cons}

With a specific focus on 2-line e.c.\ graphs, we are able to develop some constructions for producing an infinite collection of such graphs.

\begin{Theorem}
\label{cons}
Let $G$ be a 2-line e.c.\ graph. Join to $G$ a new vertex $x$ such that $x$ is adjacent to each vertex of $G$. The resulting graph $G\vee x$ is a 2-line e.c.\ graph.
\end{Theorem}

\begin{Proof}
We must verify that each pair of edges of $G\vee x$ satisfies the 2-line e.c.\ property. 
Any two edges of $G$ retain the 2-line e.c.\ property. 
For any two edges incident with $x$, say $e=\{u,x\}$ and $f=\{v,x\}$, any third edge incident with $x$ is adjacent to both $e$ and $f$, and any matching of size at least three in $G$ contains at least one edge which is adjacent to neither $e$ nor $f$. 

Now let $e=\{u,x\}$ and $f$ be an edge of $G$. 
There is an edge adjacent to both since $x$ is adjacent to each vertex of $G$ and there is an edge adjacent to neither since $f$ is a member of a matching in $G$ of size at least three. 
\end{Proof}

We can extend this result further to allow the addition of any number of new vertices to a 2-line e.c\ graph. 

\begin{Theorem}
\label{cons2}
Let $G$ be a 2-line e.c.\ graph. Join to $G$ a set $S$ of independent vertices of size $|S|\geq 2$ such that each vertex $x\in S$ is adjacent to every vertex of $G$. The resulting graph $G'$ is a 2-line e.c.\ graph.
\end{Theorem}

\begin{Proof}
We must verify that each pair of edges of $G'$ satisfies the 2-line e.c.\ property. 
By the proof of Theorem \ref{cons}, the only pairs of edges left to check are pairs of the form $e=\{x_1,u\},f=\{x_2,v\}$ where $x_1,x_2\in S$ and $u,v\in V(G)$. 

To find an edge that is adjacent to neither $e$ nor $f$, simply identify a matching of size three in $G$ and observe that $e$ and $f$ can be adjacent to at most two of the edges of the matching. 
To find an edge that is adjacent to both, take the edge $\{x_1,v\}$ if $u\neq v$ and any edge incident with $u$ in $G$ if $u=v$. 
\end{Proof}

As an immediate corollary of Theorem \ref{cons2} and by our previous observation that $K_{m,n}$ is 2-line e.c.\ for $m,n \geq 3$, we can establish the following.

\begin{Corollary}\label{multi}
Any complete multipartite graph with minimum part size at least three is a 2-line e.c.\ graph.
\end{Corollary}

We can generalise this even further to show that the join of two 2-line e.c.~graphs is itself a 2-line e.c.~graph.

\begin{Theorem}\label{join}
Let $G_1$ and $G_2$ be two 2-line e.c.~graphs each with at least three vertices.
Join $G_1$ to $G_2$ by making every vertex in $G_1$ adjacent to every vertex in $G_2$.
Then the resulting graph $G'=G_1 \vee G_2$ is a 2-line e.c.~graph.
\end{Theorem}

\begin{Proof}
We must verify that each pair of edges of $G'$ satisfies the 2-line e.c.~property.
Any two edges of $G_1$ preserve the 2-line e.c.~property, likewise for $G_2$.
By Corollary \ref{multi}, any two edges of the join have the 2-line e.c.~property as well.

Let $e=\{u,v\}$ be an edge of $G_1$ and $f=\{w,x\}$ be an edge with $w\in V(G_1)$ and $x\in V(G_2)$.
Either $\{u,x\}$ or $\{v,x\}$ serves as an edge which is adjacent to both $e$ and $f$.
Since $e$ is part of a matching of size three in $G_1$, we can find at least one edge which is adjacent to neither $e$ nor $f$.
Similar arguments verify that the 2-line e.c.~property holds between an edge of $G_2$ and an edge between $G_1$ and $G_2$.

Now let $e=\{u,v\}$ be an edge of $G_1$ and let $f=\{x,y\}$ be an edge of $G_2$.
Let $w\in V(G_1)$ and $z\in V(G_2)$ be additional distinct vertices.
Note that $\{u,x\}$ is adjacent to both $e$ and $f$, and $\{w,z\}$ is adjacent to neither $e$ nor $f$.
\end{Proof}

\section{2-Line e.c.\ Planar Graphs}\label{Planar}

Note that for a graph to be 2-line e.c., it necessarily has diameter at most 3 since any pair of disjoint edges must have a common neighbouring edge.
Small diameter in the context of planar graphs imposes upper bounds on the graph orders in terms of the maximum degree $\Delta$.

\begin{Theorem}\cite{HS}
If $G$ is planar and has diameter 2, then $|V(G)| \leq \lfloor\frac{3\Delta}{2} \rfloor + 1$ when $\Delta \geq 8$.
\end{Theorem}

\begin{Theorem}\cite{FHS}
If $G$ is planar and has diameter 3, then $|V(G)| \leq 8\Delta + 12$.
\end{Theorem}

At the same time, planar graphs cannot have high average degree (a classical application of Euler's Formula shows that the average degree must be below six) and hence they must be relatively sparse. 
In terms of diameter, the size of a planar graph is also bounded.

\begin{Theorem}\cite{FMP}
If $G$ is a connected planar graph then $|E(G)| \leq 4|V(G)| -4-3D$, where $D$ denotes the diameter of $G$.
\end{Theorem}

Despite these restrictions, infinite families of planar graphs with small diameter are known to exist, inspiring us to ask whether families of planar 2-line e.c.\ graphs exist as well.
However, as we shall see shortly, there are only finitely many such graphs as Theorem \ref{bound} establishes an upper bound on the order of a planar 2-line e.c.\ graph.
The proof of Theorem \ref{bound} relies on the following well known result of Wagner.

\begin{Theorem}\cite{Wag}
A graph $G$ is planar if and only if $G$ contains neither $K_5$ nor $K_{3,3}$ as a graph minor.
\end{Theorem}

\begin{Theorem}\label{bound}
If $G$ is a 2-line e.c.\ planar graph then $|V(G)|\leq 12$.
\end{Theorem}

\begin{Proof}
We proceed by examining the possible sizes of matchings in such a graph $G$. 
It will be useful to recall from Lemma \ref{deg3} that in a 2-line e.c.\ graph, the minimum degree, $\delta$, is at least three.
Also, since $G$ is 2-line e.c., the size of a maximum matching is at least three.

Let $M$ be a maximum matching in $G$. 
Since $G$ is 2-line e.c., each pair of edges of $M$ must have a common neighbouring edge. 
So by contracting each edge of $M$, we can observe that $G$ contains $K_{|M|}$ as a minor. 
Now if $|M|\geq 5$, then $G$ contains $K_5$ as a minor and is therefore not a planar graph. Hence $|M|\leq 4$.

Now since $M$ is maximum, every edge of $E(G)\setminus M$ must share at least one end vertex with an edge of $M$. 
Let $V_M$ be the set of end vertices of the edges of $M$. 

Now consider the set of vertices $V(G)\setminus V_M$.
We will partition this set into two disjoint sets of vertices, called Type 1 and Type 2 vertices.
Precisely, Type 1 vertices are vertices of $V(G)\setminus V_M$ which are each adjacent to both end vertices of at least one edge of $M$ and Type 2 vertices are vertices of $V(G)\setminus V_M$ which are adjacent to at most one end vertex of each edge of $M$.
Note that since $M$ is maximum, every neighbour of each vertex of $V(G)\setminus V_M$ is a member of $V_M$.

In order to prove that $|V(G)|\leq 12$, we count the total number of possible Type 1 and Type 2 vertices that $G$ could contain.
First, suppose $M=\{e_1,e_2,e_3,e_4\}$, and $u$ is a Type 1 vertex.
Without loss of generality, $G$ contains a subgraph with a structure represented by Figure \ref{type1}.

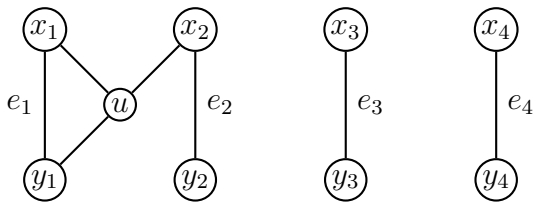
\begin{figure}[h]
\centering
\vfill
\begin{tikzpicture}[thick, main/.style={circle, draw, fill=white!75,inner sep=1pt, minimum width=12pt}]
\node[main] (v1) at (-2,2)[] {$x_1$};
\node[main] (v2) at (-2,0)[] {$y_1$};
\node[main] (v3) at (-1,1)[] {$u$};
\node[main] (v4) at (0,2)[] {$x_2$};
\node[main] (v5) at (0,0)[] {$y_2$};
\node[main] (v6) at (2,2)[] {$x_3$};
\node[main] (v7) at (2,0)[] {$y_3$};
\node[main] (v8) at (4,2)[] {$x_4$};
\node[main] (v9) at (4,0)[] {$y_4$};

\path (v1) edge node[left] {$e_1$} (v2);
\path (v1) edge node[left] {} (v3);
\path (v2) edge node[left] {} (v3);
\path (v3) edge node[left] {} (v4);
\path (v4) edge node[right] {$e_2$} (v5);
\path (v6) edge node[right] {$e_3$} (v7);
\path (v8) edge node[right] {$e_4$} (v9);

\end{tikzpicture}
\caption{The matching $M$ and a Type 1 vertex $u$.}\label{type1}
\end{figure}
Now suppose $v$ is a second Type 1 vertex.
If $v$ is adjacent to both end vertices of $e_1$, then we may form a larger matching $M^*=\{\{v,x_1\},\{y_1,u\},e_2,e_3,e_4\}$.
If $v$ is adjacent to both end vertices of $e_2$, then we may form a larger matching $M^*=\{e_1,\{u,x_2\},\{y_2,v\},e_3,e_4\}$.
Now if $v$ is adjacent to both end vertices of $e_3$ (or $e_4$), then we first observe that there must be an edge adjacent to both $e_1$ and $e_3$ (or $e_4$).
Without loss of generality, let this edge be $\{x_1,x_3\}$ (or $\{x_1,x_4\}$).
Then we may form a larger matching $M^*=\{\{y_1,u\},\{x_1,x_3\},e_2,\{y_3,v\},e_4\}$ (or $M^*=\{\{y_1,u\},\{x_1,x_4\},e_2,e_3,\{y_4,v\}\}$).
So when $|M|=4$, there can be at most one Type 1 vertex.
A similar argument shows that when $|M|=3$, there is at most one Type 1 vertex as well.

Now suppose for an edge $e=\{x,y\}$ of $M$, $G$ has two Type 2 vertices $u$ and $v$ such that $\{u,x\}$ and $\{v,y\}$ are edges. 
Then we can augment $M$ by replacing $e$ with the pair of edges $\{u,x\}$ and $\{v,y\}$.
So no two Type 2 vertices can be adjacent to both end vertices of any edge of $M$.
Consequently, if $|M|=3$, $G$ can have at most two Type 2 vertices as three or more Type 2 vertices would force a $K_{3,3}$ subgraph between the vertices of $V_M$ and the Type 2 vertices.
So if $|M|=3$, $G$ may contain at most nine vertices: six vertices of $V_M$, at most one Type 1 vertex, and at most two Type 2 vertices. 

Suppose $|M|=4$ and $u$ is a Type 1 vertex; the general structure can be observed in Figure \ref{type1}.
If $v$ is a Type 2 vertex adjacent to an end vertex of $e_1$ (say $x_1$) then we may form a larger matching $M^*=\{\{x_1,v\},\{y_1,u\}, e_2, e_3, e_4\}$.
So no Type 2 vertex can be adjacent to any end vertex of $e_1$.
Therefore, any Type 2 vertices must only be adjacent to end vertices of the edges $e_2$, $e_3$, $e_4$.
Since no two Type 2 vertices can be adjacent to both end vertices of any edge of $M$, there are precisely three vertices which are possible neighbours for a Type 2 vertex, one for each edge of $M$ other than $e_1$.
In this case, there can be at most two Type 2 vertices, since three or more would force a $K_{3,3}$ minor between the vertices of $V_M$ and the Type 2 vertices.
So if $|M|=4$ and $G$ contains a Type 1 vertex, $G$ may contain at most eleven vertices: eight vertices of $V_M$, one Type 1 vertex, and at most two Type 2 vertices.

Now suppose $|M|=4$ and $G$ contains no Type 1 vertex.
By our previous observation, no two Type 2 vertices can be adjacent to both end vertices of any edge of $M$, so a Type 2 vertex has only four possible neighbours, one end vertex from each edge of $M$.
Note that at most two Type 2 vertices may share three common neighbours among the vertices of $V_M$ since otherwise $G$ would contain a $K_{3,3}$ subgraph between the Type 2 vertices and the three common neighbours of $V_M$.
Now suppose that there are two such Type 2 vertices having three common neighbours, say $x,y,z\in V_M$.
Then any additional Type 2 vertex must share exactly two of these three common neighbours, say $x$ and $y$.
This vertex's third neighbour lies on a path with $z$ of length three consisting of two edges of $M$ and an adjacent edge shared between them.
By contracting this path, we observe a $K_{3,3}$ minor between the Type 2 vertices and the vertices of $V_M$.
Finally, assuming no pair of Type 2 vertices have three common neighbours, we conclude that there can be at most four Type 2 vertices since this is the number of distinct 3-subsets of a 4-set.
So if $|M|=4$ and $G$ contains no Type 1 vertices, $G$ may contain at most twelve vertices: eight vertices of $V_M$ and at most four Type 2 vertices.
\end{Proof}

By computer search of all planar graphs up to order 12, we established that there are precisely five 2-line e.c.\ planar graphs.
These are the graphs named Tc20, Tc30, Tc39, Tc43, and Tc44, on page 246 of \cite{Atlas}. 
For a planar representation of each graph, see \Cref{fig:a,fig:b,fig:c,fig:d,fig:e}.
This search was aided by the program $\tt{plantri}$, authored by Brinkmann and McKay \cite{plantri}.
The graphs Tc20, Tc30, and Tc39 have faces of sizes three and four whereas the graphs Tc43 and Tc44 are triangulations.
Tc43 and Tc44 are also named the heptahedral graph 34 and the Johnson solid skeleton 13 respectively.
The heptahedral graphs were first enumerated by Kirkman \cite{kirk} and Hermes \cite{herm} and the Johnson solid skeleton is the planar embedding of the pentagonal bipyramid $J_{1,3}$ \cite{johnson}.

\begin{figure}[h!]
\centering
\begin{minipage}[b][5cm][s]{.31\textwidth}
\centering
\vfill
\begin{tikzpicture}
\tikzstyle{every node}=[circle, draw, fill=white!75,inner sep=0pt, minimum width=7pt];
\node (v1) at (0,2)[] {};
\node (v2) at (-2,-1.33)[] {};
\node (v3) at (2,-1.33)[] {};
\node (v4) at (0,0.66)[] {};
\node (v5) at (-0.66,0)[] {};
\node (v6) at (0.66,0)[] {};
\node (v7) at (0,-0.66)[] {};

\draw (v1)--(v2);
\draw (v1)--(v3);
\draw (v1)--(v4);
\draw (v2)--(v3);
\draw (v2)--(v5);
\draw (v2)--(v7);
\draw (v3)--(v6);
\draw (v3)--(v7);
\draw (v4)--(v5);
\draw (v4)--(v6);
\draw (v5)--(v7);
\draw (v6)--(v7);
\end{tikzpicture}
\caption{Tc20}\label{fig:a}
\vfill
\end{minipage}\quad
\begin{minipage}[b][5cm][s]{.31\textwidth}
\centering
\vfill
\begin{tikzpicture}
\tikzstyle{every node}=[circle, draw, fill=white!75,inner sep=0pt, minimum width=7pt];
\node (v1) at (0,2)[] {};
\node (v2) at (-2,-1.33)[] {};
\node (v3) at (2,-1.33)[] {};
\node (v4) at (0,0.66)[] {};
\node (v5) at (-0.66,0)[] {};
\node (v6) at (0.66,0)[] {};
\node (v7) at (0,-0.66)[] {};

\draw (v1)--(v2);
\draw (v1)--(v3);
\draw (v1)--(v4);
\draw (v2)--(v3);
\draw (v2)--(v5);
\draw (v2)--(v7);
\draw (v3)--(v6);
\draw (v3)--(v7);
\draw (v4)--(v5);
\draw (v4)--(v6);
\draw (v5)--(v6);
\draw (v5)--(v7);
\draw (v6)--(v7);
\end{tikzpicture}
\caption{Tc30}\label{fig:b}
\vfill
\end{minipage}\quad
\begin{minipage}[b][5cm][s]{.31\textwidth}
\centering
\vfill
\begin{tikzpicture}
\tikzstyle{every node}=[circle, draw, fill=white!75,inner sep=0pt, minimum width=7pt];
\node (v1) at (0,2)[] {};
\node (v2) at (-2,-1.33)[] {};
\node (v3) at (2,-1.33)[] {};
\node (v4) at (0,0.66)[] {};
\node (v5) at (-0.66,0)[] {};
\node (v6) at (0.66,0)[] {};
\node (v7) at (0,-0.66)[] {};

\draw (v1)--(v2);
\draw (v1)--(v3);
\draw (v1)--(v4);
\draw (v1)--(v5);
\draw (v1)--(v6);
\draw (v2)--(v3);
\draw (v2)--(v5);
\draw (v2)--(v7);
\draw (v3)--(v6);
\draw (v3)--(v7);
\draw (v4)--(v5);
\draw (v4)--(v6);
\draw (v5)--(v7);
\draw (v6)--(v7);
\end{tikzpicture}
\caption{Tc39}\label{fig:c}
\vfill
\end{minipage}\quad
\begin{minipage}[b][5cm][s]{.31\textwidth}
\centering
\vfill
\begin{tikzpicture}
\tikzstyle{every node}=[circle, draw, fill=white!75,inner sep=0pt, minimum width=7pt];
\node (v1) at (0,2)[] {};
\node (v2) at (-2,-1.33)[] {};
\node (v3) at (2,-1.33)[] {};
\node (v4) at (0,0.66)[] {};
\node (v5) at (-0.66,0)[] {};
\node (v6) at (0.66,0)[] {};
\node (v7) at (0,-0.66)[] {};

\draw (v1)--(v2);
\draw (v1)--(v3);
\draw (v1)--(v4);
\draw (v1)--(v5);
\draw (v1)--(v6);
\draw (v2)--(v3);
\draw (v2)--(v5);
\draw (v2)--(v7);
\draw (v3)--(v6);
\draw (v3)--(v7);
\draw (v4)--(v5);
\draw (v4)--(v6);
\draw (v5)--(v7);
\draw (v6)--(v7);
\draw (v5)--(v6);
\end{tikzpicture}
\caption{Tc43}\label{fig:d}
\vfill
\end{minipage}\quad
\begin{minipage}[b][5cm][s]{.31\textwidth}
\centering
\vfill
\begin{tikzpicture}
\tikzstyle{every node}=[circle, draw, fill=white!75,inner sep=0pt, minimum width=7pt];
\node (v1) at (0,2)[] {};
\node (v2) at (-2,-1.33)[] {};
\node (v3) at (2,-1.33)[] {};
\node (v4) at (0,0.66)[] {};
\node (v5) at (-0.66,0)[] {};
\node (v6) at (0.66,0)[] {};
\node (v7) at (0,-0.66)[] {};

\draw (v1)--(v2);
\draw (v1)--(v3);
\draw (v1)--(v4);
\draw (v1)--(v5);
\draw (v1)--(v6);
\draw (v2)--(v3);
\draw (v2)--(v5);
\draw (v2)--(v7);
\draw (v3)--(v6);
\draw (v3)--(v7);
\draw (v4)--(v5);
\draw (v4)--(v6);
\draw (v5)--(v7);
\draw (v6)--(v7);
\draw (v4)--(v7);
\end{tikzpicture}
\caption{Tc44}\label{fig:e}
\vfill
\end{minipage}
\end{figure}

\pagebreak
\section{$n$-Line e.c.\ Hypergraphs}\label{hyper}

We know from Theorem \ref{n=2} that a graph cannot be $n$-line e.c.\ for $n\geq 3$. In order to find examples of $n$-e.c.\ line graphs for $n\geq 3$ we instead consider the line graphs of hypergraphs. 
In fact, the idea of existential closure properties in the line graphs of hypergraphs already has history in the literature, although under a more specific set of parameters. 
In particular, the block-intersection graphs of designs (which may be viewed as the line graphs of certain hypergraphs), have been studied with the $n$-existential closure property in mind. 
In a 2005 paper by Forbes, Grannell and Griggs \cite{FGG}, the block-intersection graphs of Steiner triple systems are investigated and in a 2007 paper by McKay and Pike \cite{McP}, the block-intersection graphs of more general balanced incomplete block designs were considered. 
Existential closure was also examined in the block intersection graphs of infinite designs in \cite{HPS} and \cite{PS}.

A {\it hypergraph} $H$ is a pair $(V,E)$ such that $V$ is a set of distinct vertices and $E$ is a collection of subsets of $V$ called hyperedges or simply edges. 
A hypergraph in which all edges have the same cardinality $k$ is called a {\it $k$-uniform hypergraph}.
The {\it line graph} of a hypergraph $H$, denoted $L(H)$, is the graph with vertex set $E(H)$ such that adjacency of vertices in $L(H)$ corresponds with adjacency of edges in $H$, where two edges in $H$ are adjacent if and only if they share at least one vertex. 
A matching in a hypergraph is simply a set of edges in which no two edges contain a common vertex and an independent set of vertices in a hypergraph is a set in which no two vertices are contained in a common edge.

\begin{Theorem}
Let $H$ be a hypergraph with edges of size at most $k$. 
If $H$ is $n$-line e.c.\ then $n\leq k$.
Moreover, if $H$ is a $k$-uniform $n$-line e.c.\ hypergraph, then $n\leq k$.
\end{Theorem}

\begin{Proof}
Suppose $H$ is $(k+1)$-line e.c.\ and let $e$ be an edge of $H$. By the same method as detailed in the proof of Theorem \ref{n=2}, we can build a matching of size $k+1$ in $H$ which contains $e$.

Now, since $H$ is $(k+1)$-line e.c., there must exist a $(k+2)^{\text{nd}}$ edge, adjacent to each of the previous $k+1$ edges, no two of which are adjacent. 
This is impossible, so $H$ cannot be $(k+1)$-line e.c.
Finally, in particular, this result holds when $H$ is a $k$-uniform $n$-line e.c.\ hypergraph as well.
\end{Proof}

Note that we can use existing examples of $n$-e.c.~graphs to construct $n$-line e.c.~hypergraphs for any $n$ as follows.
Let $G$ be an $n$-e.c.~graph and form a set of size $\deg(v)$ for each vertex $v\in V(G)$ consisting of the edges with which it is incident. 
From this we can build a hypergraph $H$ with $V(H)=E(G)$ and $E(H)$ consisting of the sets we have just formed.
Note that the line graph of $H$ is isomorphic to $G$, so $H$ is an $n$-line e.c.~hypergraph.
Also, if $G$ happens to be $k$-regular, then the resulting hypergraph $H$ would be $k$-uniform.

We can apply this construction to any given set of $n$-e.c.~graphs to produce additional examples of $n$-line e.c.~hypergraphs.
One explicit family of $n$-e.c.~graphs is the set of Paley graphs.
In \cite{BEH} and \cite{BT} it was shown that for any $n$, every sufficiently large Paley graph is $n$-e.c.
Since Paley graphs are necessarily regular, we conclude that there exist sufficiently large uniform $n$-line e.c.~hypergraphs for any $n$.

If we once again focus on $n=2$, then we can find analogous results for hypergraphs to those presented in Section~\ref{Cons}.

\begin{Theorem}\label{XY}
Let $X$ and $Y$ be disjoint sets of vertices. 
Let $H$ be the hypergraph on $X\cup Y$ along with all possible edges of size $k$ such that each edge has a non-empty intersection with both $X$ and $Y$.
If $|X|\geq|Y|\geq 2k-1$, then $H$ is a $k$-uniform 2-line e.c.~hypergraph.
\end{Theorem}

\begin{Proof}
We must verify that each pair of edges in $H$ satisfies the 2-line e.c.~property.
Let $e=\{u_1,u_2,\dots, u_k\}$ and $f=\{v_1,v_2,\dots, v_k\}$ be two edges of $H$.
Without loss of generality we may assume that $u_1,v_1\in X$ and $u_k,v_k\in Y$.

Any third edge which contains $u_1$ and $v_1$ is adjacent to both $e$ and $f$.
Any edge $\{x_1,x_2,\dots,x_k\}$ where each $x_i$ is distinct from each $u_i$ and each $v_i$, is adjacent to neither $e$ nor $f$.
This is possible since $|X|\geq|Y|\geq 2(k-1)+1=2k-1$.
Finally, since $e\neq f$, there exists at least one vertex of $e$ which is not contained in $f$ and vice versa.
So pick this vertex and observe that at least one of its incident edges will serve as an edge which is adjacent to $e$ and not $f$ and vice versa.
\end{Proof}

Using Theorem~\ref{XY} and similar arguments used in the proof of Theorem~\ref{join} we can establish the following corollary.

\begin{Corollary}
Let $H_1$ and $H_2$ be $k$-uniform 2-line e.c.~hypergraphs on distinct sets of vertices. 
Let $H$ be the hypergraph $H_1\cup H_2$ along with all possible edges of size $k$ such that each edge has a non-empty intersection with both $V(H_1)$ and $V(H_2)$.
If $|V(H_1)|\geq|V(H_2)|\geq 2k-1$, then $H$ is a $k$-uniform 2-line e.c.~hypergraph.
\end{Corollary}

\section{Acknowledgements}

Authors Burgess and Pike acknowledge NSERC Discovery Grant support and Luther acknowledges NSERC scholarship support.



\end{document}